\documentclass[a4paper,12pt]{amsart}

\usepackage{amsmath}
\usepackage{amssymb}
\usepackage{mathrsfs}
\usepackage{ifthen}
\usepackage{graphicx}
\usepackage[T1]{fontenc} 
\usepackage{enumerate}

\def\switchlinenumbers{\@ifstar
    {\let\makeLineNumberOdd\makeLineNumberRight
     \let\makeLineNumberEven\makeLineNumberLeft}%
    {\let\makeLineNumberOdd\makeLineNumberLeft
     \let\makeLineNumberEven\makeLineNumberRight}%
    }

\def\setmakelinenumbers#1{\@ifstar
  {\let\makeLineNumberRunning#1%
   \let\makeLineNumberOdd#1%
   \let\makeLineNumberEven#1}%
  {\ifx\c@linenumber\c@runninglinenumber
      \let\makeLineNumberRunning#1%
   \else
      \let\makeLineNumberOdd#1%
      \let\makeLineNumberEven#1%
   \fi}%
  }

\nonstopmode \numberwithin{equation}{section}
\newtheorem*{theorem*}{Theorem}

\newtheorem{thm}{Theorem}[section]
\newtheorem{cor}[equation]{Corollary}
\newtheorem{lem}[equation]{Lemma}

\theoremstyle{definition}
\newtheorem{defn}{Definition}[section]

\newtheorem{qsn}[equation]{Question}
\newtheorem{prob}[equation]{Problem}

\newtheorem{innercustomthm}{Theorem}
\newenvironment{customthm}[1]
  {\renewcommand\theinnercustomthm{#1}\innercustomthm}
  {\endinnercustomthm}

\newcounter{minutes}
\divide\time by 60
\newcounter{hours}
\multiply\time by 60
\addtocounter{minutes}{-\time}

\newcounter {own}
\def\theown {\thesection       .\arabic{own}}

\newenvironment{pf}[1][]{%
 \vskip 3mm
 \noindent
 \ifthenelse{\equal{#1}{}}%
  {{\slshape Proof. }}%
  {{\slshape #1.} }%
 }%
{\qed\bigskip}

\newcounter{alphabet}

\def\be{\begin{equation}}
\def\ee{\end{equation}}

\newcommand{\bee}{\begin{enumerate}}
\newcommand{\eee}{\end{enumerate}}

\newcommand{\blem}{\begin{lem}}
\newcommand{\elem}{\end{lem}}
\newcommand{\bthm}{\begin{thm}}
\newcommand{\ethm}{\end{thm}}
\newcommand{\bcor}{\begin{cor}}
\newcommand{\ecor}{\end{cor}}
\newcommand{\beg}{\begin{examp}}
\newcommand{\eeg}{\end{examp}}
\newcommand{\begs}{\begin{examples}}
\newcommand{\eegs}{\end{examples}}
\newcommand{\bdefe}{\begin{defin}}
\newcommand{\edefe}{\end{defin}}
\newcommand{\bprob}{\begin{prob}}
\newcommand{\eprob}{\end{prob}}
\newcommand{\bei}{\begin{itemize}}
\newcommand{\eei}{\end{itemize}}
\newcommand{\real}{{\operatorname{Re}\,}}

\newcommand{\norm}[1]{\left\lVert#1\right\rVert}
\newcommand{\abs}[1]{\left\lvert#1\right\rvert}
\newcommand{\innpdct}[1]{\left\langle#1\right\rangle}

\begin{document}

\title{Operator valued analogues of multidimensional Bohr's inequality}

\author{Vasudevarao Allu}
\address{Vasudevarao Allu,
School of Basic Sciences,
Indian Institute of Technology Bhubaneswar,
Bhubaneswar-752050, Odisha, India.}
\email{avrao@iitbbs.ac.in}

\author{Himadri Halder}
\address{Himadri Halder,
School of Basic Sciences,
Indian Institute of Technology Bhubaneswar,
Bhubaneswar-752050, Odisha, India.}
\email{himadrihalder119@gmail.com}

\subjclass[{AMS} Subject Classification:]{Primary 32A05, 32A10,47A56,; Secondary 47A63,30H05}

\def\thefootnote{}
\footnotetext{ {\tiny File:~\jobname.tex,
printed: \number\year-\number\month-\number\day,
          \thehours.\ifnum\theminutes<10{0}\fi\theminutes }
} \makeatletter\def\thefootnote{\@arabic\c@footnote}\makeatother

\begin{abstract}
Let $\mathcal{B}(\mathcal{H})$ be the algebra of all bounded linear operators on a complex Hilbert space $\mathcal{H}$.
In this paper, we first establish several sharp improved and refined versions of the Bohr's inequality for the functions in the class $H^{\infty}(\mathbb{D},\mathcal{B}(\mathcal{H}))$ of bounded analytic functions from the unit disk $\mathbb{D}:=\{z \in \mathbb{C}:|z|<1\}$ into $\mathcal{B}(\mathcal{H})$. For the complete circular domain $Q \subset \mathbb{C}^n$, we prove the multidimensional analogues of the operator valued Bohr's inequality established by G. Popescu [Adv. Math. 347 (2019), 1002-1053]. Finally, we establish the multidimensional analogues of several improved Bohr's inequalities for operator valued functions in $Q$.	
\end{abstract}

\maketitle
\pagestyle{myheadings}
\markboth{Vasudevarao Allu and  Himadri Halder}{Operator valued analogues of multidimensional Bohr's inequality}

\section{Introduction and some basic questions}
Let $H^{\infty}(\mathbb{D},X)$ be the space of bounded analytic functions from the unit disk $\mathbb{D}:=\{z \in \mathbb{C}:|z|<1\}$ into a complex Banach space $X$ with $\norm{f}_{H^{\infty}(\mathbb{D},X)}:=\sup_{|z|<1} \norm{f(z)}$. Let $\mathcal{B}(\mathbb{D},X)$ be the class of functions $f$ in $H^{\infty}(\mathbb{D},X)$ with $\norm{f}_{H^{\infty}(\mathbb{D},X)} \leq 1$. The Bohr radius $R(X)$ for the class $\mathcal{B}(\mathbb{D},X)$ is defined by (see \cite{Blasco-OTAA-2010})
\begin{equation*}
R(X):=\sup \left\{r \in (0,1): M_{r}(f)\leq 1 \,\, \mbox{for all} \,\, f(z)= \sum_{k=0}^{\infty} x_{k} z^k \in \mathcal{B}(\mathbb{D},X), z \in \mathbb{D}\right\},
\end{equation*}
where $M_{r}(f)=\sum_{k=0}^{\infty} \norm{x_{k}} r^k$ is the associated majorant series of $f \in H^{\infty}(\mathbb{D},X)$. 
The remarkable theorem of Harald Bohr \cite{Bohr-1914} (in improved form) states that $R(X)=1/3$ for $X=\mathbb{C}$, where the norm of $X$ is the usual modulus of complex numbers. The interest in the Bohr's theorem has been revived when Dixon \cite{Dixon & BLMS & 1995} used it to answer a long-standing question on the characterization of the Banach algebras satisfying the non-unital von Neumann inequality. For the last two decades, there has been an extensive research carried out to the extensions of analytic functions of several complex variables, to planar harmonic mappings, to polynomials, to solutions of elliptic partial differential equations, and to more abstract settings. In $1997$, Boas and Khavinson \cite{boas-1997} introduced the $n$-dimensional Bohr radius $K_{n}$ for the Hardy space of bounded analytic functions on the unit polydisk, and obtained the upper and lower bounds of $K_{n}$. In $2006$, an improved version of the lower estimate of $K_{n}$ was obtained by Defant and Frerick \cite{defant-2006}. Further estimation of $K_{n}$ has been obtained by Defant {\it et al.} \cite{defant-2011} by using the hypercontractivity of the polynomial Bohnenblust-Hille inequality. In $2014$, Bayart {\it et al.} \cite{bayart-advance-2014} obtained the exact asymptotic behaviour of $K_{n}$. In $2019$, Popescu \cite{popescu-2019} extended the Bohr inequality for free holomorphic functions to polyballs. For more interesting aspects and generalization of multidimensional Bohr's inequality, we refer to \cite{aizenberg-2001,aizn-2007,alkhaleefah-2019,Ayt & Dja & BLMS & 2013,Djakov & Ramanujan & J. Anal & 2000,Liu-Pon-PAMS-2020}. In $2004$, Paulsen and Singh \cite{paulsen-2004-PAMS} extended Bohr's theorem to Banach algebras by finding a general version of Bohr inequality which is valid in the context of uniform algebras. In $2021$, Bhowmik and Das \cite{bhowmik-2021} extensively studied Bohr inequality for operator valued functions. Further results on Bohr radius, we refer to \cite{Himadri-Vasu-P1,Himadri-Vasu-P2,bene-2004,Bohr-Blaschke-Book}.
\vspace{1mm}

For $p \in [1, \infty)$, let $H^{p}(\mathbb{D},X)$ be the space of analytic functions from $\mathbb{D}$ into a complex Banach space $X$ such that 
\begin{equation} \label{him-p7-e-1.11}
\norm{f}_{H^{p}(\mathbb{D},X)}= \sup \limits _{0<r<1} \left(\int_{0}^{2\pi} \norm{f(re^{it})}^p \frac{dt}{2\pi}\right)^{1/p} < \infty.
\end{equation}
In \cite{bene-2004}, B\'{e}n\'{e}teau {\it et al.} have shown that there is no Bohr phenomenon in the Hardy spaces $\norm{f}_{H^{p}(\mathbb{D},\mathbb{C})}$ for $1 \leq p <\infty$. In fact, they have shown that there is no Bohr phenomenon in complex valued Hardy spaces $H^q$ for $q \in (0,\infty)$. In \cite{Djakov & Ramanujan & J. Anal & 2000}, Djakov and Ramanujan have extensively studied the $p$-Bohr phenomenon for the power series of the form $\sum_{k=0}^{\infty} |a_{k}|^p r^k$ for $p \in [1,\infty)$, where $f(z)=\sum_{k=0}^{\infty} a_{k}z^k$ is a bounded analytic function in $\mathbb{D}$. In \cite{Djakov & Ramanujan & J. Anal & 2000}, the notion of $p$-Bohr inequality has been extended to the bounded analytic functions of several variables. 
\par
In this paper, we establish multidimensional analogues of several improved $p$-Bohr's inequalities for the operator valued bounded analytic functions.
\subsection{Bohr theorem for operator valued analytic function} 
One of the main aims of the present paper is to study Bohr inequality in the setting of operator valued analytic functions in the unit disk $\mathbb{D}$, to be more specific, for functions in $H^{\infty}(\mathbb{D},X)$, where $X=\mathcal{B}(\mathcal{H})$ is the algebra of all bounded linear operators on a complex Hilbert space $\mathcal{H}$. For the rest of our discussion on this, we need to fix some basic notations. For $T \in \mathcal{B}(\mathcal{H})$, $\norm{T}$ denotes the operator norm of $T$. The adjoint operator $T^{*}:\mathcal{H} \rightarrow \mathcal{H}$ of $T$ is defined by $\innpdct{Tx,y}=\innpdct{x,T^{*}y}$ for all $x, y \in \mathcal{H}$. The operator $T$ is said to be normal if $T^{*}T=TT^{*}$, self-adjoint if $T^{*}=T$, and positive if $\innpdct{Tx,x} \geq 0$ for all $x \in \mathcal{H}$. The absolute value of $T$ is defined by $\abs{T}:=\left(T^{*}T\right)^{1/2}$, while $S^{1/2}$ denotes the unique positive square root of a positive operator $S$. 
Let $I$ be the identity operator on $\mathcal{H}$. Let $f \in H^{\infty}(\mathbb{D},\mathcal{B}(\mathcal{H}))$ be a bounded analytic function with the expansion
\begin{equation} \label{P9-e-1.1}
f(z)= \sum_{n=0}^{\infty} A_{n}z^n \,\,\,\, \mbox{for} \,\,\,\, z \in \mathbb{D},
\end{equation}
where $A_{n} \in \mathcal{B}(\mathcal{H})$ for all $n \in \mathbb{N} \cup \{0\}$. For each $f \in H^{\infty}(\mathbb{D},\mathcal{B}(\mathcal{H}))$ of the form \eqref{P9-e-1.1}, the function $r \mapsto M_{r}(f)$ is an increasing function in $[0,1)$ with $m_{0}(f)=\norm{A_{0}} \leq 1$, where $M_{r}(f)$ is the associated majorant series of $f$ defined by $M_{r}(f)=\sum_{k=0}^{\infty} \norm{A_{k}} r^k$ for $r \in [0,1)$. For each fixed $z \in \mathbb{D}$, we denote
$
\mathcal{G}_{z}:= \left\{f(z)=\sum_{k=0}^{\infty} A_{k}z^k: f\in H^{\infty}(\mathbb{D},\mathcal{B}(\mathcal{H}))\right\}.
$
In \cite{Himadri-Vasu-P8}, Allu and Halder have proved that the space $\mathcal{G}_{z}$ with norm $M_{r}$ constitutes a Banach algebra and have shown that functions in $\mathcal{G}_{z}$ satisfy a von Neumann type inequality. Set $\chi=\mathcal{B}(\mathbb{D},\mathcal{B}(\mathcal{H}))$. By the similar definition as in \cite{Bomberi-1962} for the complex valued functions, we define
\begin{equation} \label{P9-e-1.2}
m(\chi,r):=\sup_{f\in \chi} M_{r}(f)\,\,\,\,\, \mbox{for} \,\,\,\, r\in [0,1).
\end{equation} 
Clearly, $m(\chi,0)=1$. It is worth mentioning that $m(\chi,r)$ is an increasing function of $r$ and hence $m(\chi,r)\geq 1$ for $r \in [0,1)$. For the arbitrary functions $f$ in $H_{\infty}(\mathbb{D},\mathcal{B}(\mathcal{H}))$, not necessarily $\norm{f}_{H^{\infty}(\mathbb{D},\mathcal{B}(\mathcal{H}))} \leq 1$, we have $m(\chi,r) \leq M_{r}(f) \norm{f}_{H^{\infty}(\mathbb{D},\mathcal{B}(\mathcal{H}))}$. Thus, finding or estimating $m(\chi,r)$, it is relevance to know $\norm{f(z)}$ {\it i.e.,} to understand the rate of the growth of functions in $H^{\infty}(\mathbb{D},\mathcal{B}(\mathcal{H}))$. However, a precise value for $m(\chi,r)$ is not known for all $r \in [0,1)$. But we can estimate the bounds for $m(\chi,r)$. In \cite{Bomberi-1962}, Bombieri has obtained upper bounds and lower bounds for $m(\chi ,r)$ for complex valued bounded analytic functions. \par

In Lemma \ref{P9-lem-2.1}, we obtain an operator valued analogue of Bombieri's bound of $m(\chi ,r)$. We now generalize the notion of Bohr radius for the class $z \chi=\{zf:f \in \chi\}$ {\it i.e.,} for the class of functions $f \in \chi$ with $f(0)=0$. The Bohr radius for the class $z \chi$ is the largest radius $R_{1}$ such that (a) $M_{r}(f)=\sum_{k=0}^{\infty}\norm{A_{k}}r^{k+1} \leq 1$ for $r \in [0,R_{1}]$ and for all $f(z)=\sum_{k=0}^{\infty} A_{k}z^{k+1} \in z\chi$, (b) when $r \in (R_{1},1)$, there is a function $f \in z\chi$ such that $M_{r}(f)>1$. In view of \eqref{P9-e-1.2}, it is easy to see that $R_{1}$ is the unique root of 
\begin{equation} \label{P9-e-1.3}
r\, m(\chi,r)=1.
\end{equation}
Since $m(\chi,r)$ is strictly increasing, the function $r \mapsto r\, m(\chi,r)$ is also strictly increasing in $r\in [0,1)$, which shows that \eqref{P9-e-1.3} has the unique root in $(0,1)$. Bombieri \cite{Bomberi-1962} has proved that $R_{1}=1/\sqrt{2}$ for the complex valued bounded functions in $\mathbb{D}$. Later, Paulsen {\it et al.} \cite{paulsen-2002} have extensively studied the radius $R_{1}$ for complex valued functions. In the present paper, we obtain $R_{1}$ for operator valued functions in $\chi=\mathcal{B}(\mathbb{D}, \mathcal{B}(\mathcal{H}))$ in Lemma \ref{P9-lem-2.2}.

In $2019$, Popescu \cite{popescu-2019} proved the following interesting result, which is an analogue of the classical theorem of Bohr for operator valued bounded analytic functions in $\mathbb{D}$. 
\begin{customthm}{A} \cite{popescu-2019}\label{him-vasu-P9-thm-1.1} 
	Let $f \in H^{\infty}(\mathbb{D},\mathcal{B}(\mathcal{H}))$ be an operator valued bounded holomorphic function with the expansion \eqref{P9-e-1.1} such that $A_{0}=a_{0}I$, $a_{0} \in \mathbb{C}$. Then 
	 \begin{equation} \label{him-vasu-P9-e-1.1}
	 \sum _{n=0}^{\infty} \norm{A_{n}} \, r^n \leq \norm{f}_{H^{\infty}(\mathbb{D},\mathcal{B}(\mathcal{H}))} \,\,\,\,\,\, \mbox{for} \,\,\,\,\, |z|=r \leq \frac{1}{3}
	 \end{equation}
	 and $1/3$ is the best possible constant. Moreover, the inequality is strict unless $f$ is a constant.
\end{customthm}
The proof of Theorem \ref{him-vasu-P9-thm-1.1} relies on the bound of the norm of coefficients $\norm{A_{n}}$, which may be obtained as an application of the operator counterpart of the Schwarz-Pick inequality. 
\begin{lem} \cite{anderson-2006} (counterpart of Schwarz-Pick inequality)\label{him-vasu-P6-lem-2.1}
	Let $B(z)$ be an analytic function with values in $\mathcal{B}(\mathcal{H})$ and satisfying $\norm{B(z)} \leq 1$ on $\mathbb{D}$. Then
	\begin{equation*} 
	(1-|a|)^{n-1}\, \norm{\frac{B^{(n)}(a)}{n!}} \leq \frac{\norm{I-B(a)^{*}B(a)}^{1/2} \, \norm{I-B(a)B(a)^{*}}^{1/2}}{1-|a|^2} 
	\end{equation*}
	for each $a \in \mathbb{D}$ and $n=1,2,\ldots$.
\end{lem}

We note that for $f \in \mathcal{B}(\mathbb{D},\mathcal{B}(\mathcal{H}))$ of the form \eqref{P9-e-1.1} with $A_{0}=a_{0}I$, $|a_{0}|<1$. Without loss of generality, assume that $\norm{f}_{H^{\infty}(\mathbb{D},\mathcal{B}(\mathcal{H}))} \leq 1$. Then, by the virtue of Lemma \ref{him-vasu-P6-lem-2.1}, putting $a=0$, we obtain
\begin{equation} \label{him-vasu-P9-e-1.3}
\norm{A_{n}} \leq \norm{I-|A_{0}|^2}=1-|a_{0}|^2 \,\,\,\,\, \mbox{for} \,\,\,\, n\in \mathbb{N}.
\end{equation}
In view of \eqref{him-vasu-P9-e-1.3}, the proof of Theorem \ref{him-vasu-P9-thm-1.1} follows immediately. For the sharpness of the constant $1/3$, we consider the following function 
\begin{equation} \label{P9-e-1.7}
\psi_{a}(z)= \left(\frac{a-z}{1-az}\right)I= A_{0} + \sum_{k=1}^{\infty} A_{k} z^k \,\,\, \, z \in \mathbb{D},
\end{equation}
where $A_{0}= aI, A_{k}=-(1-a^2)a^{k-1} I$ and for some $a \in [0,1)$. Then it is easy to see that $M_{r}(f) =\sum_{k=0}^{\infty} \norm{A_{k}}r^k= (a+(1-2a^2)r)/(1-ar)> 1$ whenever $r > 1/(1+2a)$. Taking $a$ is very close to $1$ {\it i.e.,} $a \rightarrow 1^{-}$, we obtain $M_{r}(f)>1$ for $r>1/3$, which shows that $1/3$ is the best possible.
\vspace{1mm}

 We say that $\mathcal{B}(\mathbb{D},\mathcal{B}(\mathcal{H}))$ satisfies Bohr phenomenon if all the functions in $\mathcal{B}(\mathbb{D},\mathcal{B}(\mathcal{H}))$ satisfies the inequality \eqref{him-vasu-P9-e-1.1} for $r \leq 1/3$. It is worth mentioning that the constant $1/3$ does not dependent on the coefficients of functions. 
\par
In the recent years, refined versions and improved versions of Bohr inequality in the case of bounded analytic functions have become central research interest in one and several complex variables. Several authors have established various refined and improved versions of Bohr's inequality (see \cite{Ahamed-Allu-Halder-P3-2020,Kay & Pon & AASFM & 2019}). We now recall some of them. Let $f \in H^{\infty}(\mathbb{D},\mathbb{C})$ be of the form 
\begin{equation} \label{P9-e-1.5}
f(z)=\sum_{k=0}^{\infty} a_{k}z^k \,\,\,\, \mbox{for} \,\,\,\, z \in \mathbb{D}. 
\end{equation}
Let $S_{r}$ be the area of the image of the subdisk $\mathbb{D}_{r}=\{z:|z|<r\}$ under the map $f$ given by \eqref{P9-e-1.5}. Then it is known that $S_{r}/\pi= \sum_{k=1}^{\infty}k|a_{k}|^2 r^{2k}$ (see \cite{Evd-Ponn-Rasi-2020}).
\begin{customthm}{B}  \label{P9-thm-B} Let $f \in \mathcal{B}(\mathbb{D},\mathbb{C})$ of the form \eqref{P9-e-1.5}. Then 
\begin{enumerate}
	\item [\rm{(a)}]  $
	\sum\limits_{k=0}^{\infty} |a_{k}|r^k + \left(\frac{1}{1+|a_{0}|} + \frac{r}{1-r}\right) \sum\limits_{k=1}^{\infty} |a_{k}|^2 r^{2k} \leq 1 
	$ for $|z|=r\leq 1/3$.
	The constant $1/3$ is the best possible.
	\item [\rm{(b)}] $\sum \limits_{k=0}^{\infty}|a_{k}|r^k+\frac{8}{9}\left(\frac{S_{r}}{\pi}\right)\leq 1\;\;\mbox{for}\;\; r\leq 1/3$. The bound $8/9$ and the constant $1/3$ cannot be replaced by a larger quantity.\\
	\item [\rm{(c)}] $
	|f(z)| + \sum \limits_{k=N}^{\infty} \abs{a_{k}}r^k \leq 1 \,\,\,\,\,\,\mbox{for} \,\,\,\,\,|z|=r \leq R_{N,1},
	$
	where $R_{N,1}$ is the unique root in $(0,1)$ of $2(1+r) r^N - (1-r)^2=0$. The constant $R_{N,1}$ is the best possible.
\end{enumerate}
\end{customthm}
It is important to note that the proof of the part (a) and part (b) can be obtained by putting $\gamma=0$ in \cite[Theorem 2]{Evd-Ponn-Rasi-2020} and \cite[Theorem 1]{Evd-Ponn-Rasi-2020} respectively. The proof of the part (c) can be obtained by putting $\gamma=0$ in \cite[Theorem 2.7]{Ahamed-Allu-Halder-P3-2020}.

Like the quantity $S_{r}/\pi= \sum_{k=1}^{\infty}k|a_{k}|^2 r^{2k}$ for the functions $f(z)=\sum_{k=0}^{\infty} a_{k} z^k$ in $ H^{\infty}(\mathbb{D},\mathbb{C})$, we define 
\begin{equation} \label{P9-e-2.5-adjust}
S(z)=\sum_{k=1}^{\infty}k \norm{A_{k}}^2 r^{2k}, \,\,\, |z|=r<1,
\end{equation}
for the functions $f(z)=\sum_{k=0}^{\infty} A_{k} z^k \in H^{\infty}(\mathbb{D},\mathcal{B}(\mathcal{H}))$. 
Let $ G(w) $ be a polynomial defined by
\begin{equation}\label{e-1.8}
G(w)=c_1w+c_2w^2+\cdots+c_lw^l\;\; \mbox{for}\;\;  c_i\in\mathbb{R}^{+},\;\; i=1, 2, \cdots, l.
\end{equation}

\subsection{Multidimensional analogues of operator valued Bohr's inequality}
Let $\alpha$ be an $n$-tuple $(\alpha_{1},\alpha_{2},\ldots,\alpha_{n})$ of nonnegative integers, $|\alpha|$ be the sum $\alpha_{1}+\cdots + \alpha_{n}$ of its components, $\alpha !$ denotes the product $\alpha_{1}! \alpha_{2}!\ldots \alpha_{n}!$, $z$ denotes an $n$-tuple $(z_{1},z_{2},\ldots, z_{n})$ of complex numbers, and $z^{\alpha}$ denotes the product $z_{1}^{\alpha_{1}} z_{2}^{\alpha_{2}}\ldots z_{n}^{\alpha_{n}}$. Using the standard multi-index notation, we write an operator valued $n$-variable power series 
\begin{equation} \label{P9-e-1.10}
f(z)=\sum_{\alpha} A_{\alpha} z^{\alpha}, \,\,\,\,\, A_{\alpha} \in \mathcal{B}(\mathcal{H}).
\end{equation}
Let $\mathbb{D}^n=\{z \in \mathbb{C}^n: z=(z_{1}, \ldots,z_{n}), |z_{j}|<1, j=1,2,\ldots,n\}$.
Let $K_{n}(\mathcal{H})$ be the largest nonnegative number such that the power series \eqref{P9-e-1.10} converges in $\mathbb{D}^n$ and $\norm{f}_{H^{\infty}(\mathbb{D},\mathcal{B}(\mathcal{H}))} \leq 1$, then 
\begin{equation}
\sum_{\alpha} \norm{A_{\alpha}} \abs{z^{\alpha}}\,\,\,\, \mbox{for all }\,\,\,\,\, z\in K_{n}(\mathcal{H}).\mathbb{D}^n.
\end{equation}
\begin{defn}
A domain $D \subset \mathbb{C}^n$ is said to be a Reinhardt domain centered at $0 \in D$ if for any $z=(z_{1},\ldots,z_{n})\in D$, and for each $\theta_{k}\in [0,2\pi]$, $k=1,2,\ldots,n$, we have that $(z_{1}e^{i\theta_{1}},\ldots,z_{n}e^{i\theta_{n}}) \in D$. We say that $D \subset \mathbb{C}^n$ is a complete Reinhardt domain if for each $z=(z_{1},\ldots,z_{n})\in D$, and for each $|\xi_{k}|\leq1$, $k=1,2,\ldots,n$, we have that $\xi .z=(z_{1}\xi_{1}, \ldots,z_{n}\xi_{n}) \in D$.
\end{defn}
A domain $Q \subset \mathbb{C}^n$ is called a circular domain centered at $0\in Q$ if for any $z=(z_{1},\ldots,z_{n})\in Q$, and for each $\theta\in [0,2\pi]$, we have that $(z_{1}e^{i\theta},\ldots,z_{n}e^{i\theta}) \in Q$. We say that $Q \subset \mathbb{C}^n$ is a complete circular domain if for each $z=(z_{1},\ldots,z_{n})\in Q$, and for each $|\xi|\leq1$, we have that $\xi .z=(z_{1}\xi_{1}, \ldots,z_{n}\xi_{n}) \in Q$. For a complete circular domain $Q\subseteq \mathbb{C}^n$ centered at $0\in Q$, every analytic function in $Q$ can be expressed into the following homogeneous polynomials 
\begin{equation} \label{him-vasu-P9-e-1.2}
f(z)= \sum_{k=0}^{\infty} P_{k}(z)\,\,\,\,\,\, \mbox{for} \,\,\,\, z \in Q,
\end{equation}
where $P_{k}(z)=\sum_{|\alpha|=k} A_{\alpha} z^{\alpha}$ is a homogeneous polynomial of degree $k$ and $P_{0}(z)=f(0)$. 
\vspace{1mm}

It is natural to raise the following question.
\begin{qsn} \label{qsn-1.11}
	Can we establish the analogue of Theorem \ref{P9-thm-B}   for the operator valued functions in $\mathcal{B}(\mathbb{D},\mathcal{B}(\mathcal{H}))$? If so, then what is the multidimensional analogue of Theorem \ref{P9-thm-B} for the functions in $\mathcal{B}(\mathbb{D},\mathcal{B}(\mathcal{H}))$?
\end{qsn}
In Theorem \ref{P9-thm-B} (b), it is worth mentioning that Bohr inequality is improved by adding one degree polynomial in $S_{r}/\pi$ with the majorant  series $\sum_{k=0}^{\infty} |a_{k}|r^k$. Therefore, we have the following question.
\begin{qsn} \label{qsn-1.12}
	Is it possible to establish an improved version of Theorem \ref{P9-thm-B} (b) by adding the polynomial $G\left(S_{r}/\pi\right)$ with the majorant  series $\sum_{k=0}^{\infty} ||A_{k}||r^k$ for $f \in \mathcal{B}(\mathbb{D},\mathcal{B}(\mathcal{H}))$? If so, then what is the multidimensional analogue of Theorem \ref{P9-thm-B} for the functions in $\mathcal{B}(\mathbb{D},\mathcal{B}(\mathcal{H}))$?
\end{qsn}
One of the main aims of this paper is to answer Question \ref{qsn-1.11} and Question \ref{qsn-1.12} affirmatively.
\section{Main Results}

\begin{thm} \label{him-vasu-P9-thm-2.1}
If the series \eqref{him-vasu-P9-e-1.2} converges in the domain $Q$ such that the estimate $\norm{f(z)} <1$ holds in $Q$ and $f(0)=a_{0}I$, $a_{0} \in \mathbb{C}$, $|a_{0}|<1$, then 
\begin{equation} \label{him-vasu-P9-e-2.1}
\sum_{k=0}^{\infty} \norm{P_{k}(z)} <1
\end{equation}
in the homothetic domain $(1/3)Q$. Moreover, if $Q$ is convex, the constant $1/3$ is the best possible.
\end{thm}
\noindent In particular, if $f$ is complex valued bounded analytic functions in the domain $Q$, we can obtain the multidimensional analogues of Bohr's inequality for complex valued bounded analytic functions in the domain $Q$, which has been independently proved by Aizenberg \cite{aizn-2000}.
We obtain the Bohr radius for the functions $f \in \mathcal{B}(\mathbb{D}, \mathcal{B}(\mathcal{H}))$ with the initial coefficients $f(0)=0$ in the following lemma. Recall that $R_{1}$ is the unique root of \eqref{P9-e-1.3}.
\begin{lem} \label{P9-lem-2.7}
	Let $f \in \chi=\mathcal{B}(\mathbb{D}, \mathcal{B}(\mathcal{H}))$ with $f(0)=0$. Then $R_{1}=1/\sqrt{2}$.
\end{lem}
 If $f \in \mathcal{B}(\mathbb{D}, \mathbb{C})$, then the well known Bombieri \cite{Bomberi-1962} result follows from Lemma \ref{P9-lem-2.7}. In the following result, we obtain an multidimensional analogues of Lemma \ref{P9-lem-2.7}. 
\begin{thm} \label{him-vasu-P9-thm-2.2}
	If the series \eqref{him-vasu-P9-e-1.2} converges in the domain $Q$ with the estimate $\norm{f(z)} <1$ holds in $Q$ and $f(0)=0$, then 
	\begin{equation} \label{him-vasu-P9-e-2.2}
	\sum_{k=1}^{\infty} \norm{P_{k}(z)} <1
	\end{equation}
	in the homothetic domain $(1/\sqrt{2})Q$. Moreover, if $Q$ is convex, the constant $1/\sqrt{2}$ is the best possible.
\end{thm}

\begin{cor}
Suppose that the series \eqref{P9-e-1.10} converges in the polydisk $\mathbb{D}^n$ such that $f(0)=0$ and $\norm{f(z)}<1$ in $\mathbb{D}^n$. Then 
\begin{equation*}
\sum_{k=1}^{\infty} \norm{\sum_{|\alpha|=k}A_{\alpha}z^{\alpha}} \leq 1
\end{equation*}
in the polydisk $(1/\sqrt{2})\mathbb{D}^n$. The constant $1/\sqrt{2}$ is the best possible. 
\end{cor}
We establish the multidimensional analogues of Theorem \ref{P9-thm-B} (c) for the operator valued analytic functions in the complete circular domain $Q$.
\begin{thm} \label{P9-thm-2.3}
Suppose that $Q$ is a complete circular domain centered at $0 \in Q \subset \mathbb{C}^{n}$. If the series \eqref{him-vasu-P9-e-1.2} converges in $Q$ such that $\norm{f(z)}<1$ for all $z \in Q$ and $f(0)=a_{0}I$, $a_{0} \in \mathbb{C}$, $|a_{0}|<1$. Then for $p \in (0,1]$, we have
\begin{equation} \label{P9-e-2.4}
\norm{f(z)}^p + \sum_{k=N}^{\infty} \norm{P_{k}(z)} \leq 1
\end{equation}
in the homothetic domain $(R_{N,p})Q$, where $R_{N,p}$ is the positive root in $(0,1)$ of the equation 
\begin{equation} \label{P9-e-2.5}
2(1+r) r^N - p(1-r)^2=0.
\end{equation}
Moreover, if $Q$ is convex, then the constant $R_{N,p}$ is the best possible.
\end{thm}
In the following result, we prove the multidimensional version of Theorem \ref{P9-thm-B} (a) for the operator valued analytic functions in the complete circular domain $Q$.
\begin{thm} \label{P9-thm-2.4}
Assume that the series \eqref{him-vasu-P9-e-1.2} converges in the domain $Q$ such that $\norm{f(z)} <1$ for all $z \in Q$ and $f(0)=a_{0}I$, $a_{0} \in \mathbb{C}$, $|a_{0}|<1$. Then 
\begin{equation} \label{P9-e-2.5-a}
\sum_{k=0}^{\infty} \norm{P_{k}(z)} + \left(\frac{1}{1+\norm{f(0)}} + \frac{r}{1-r}\right) \sum_{k=1}^{\infty} \norm{P_{k}(z)}^2 \leq 1.
\end{equation}
holds in the homothetic domain $(1/3)Q$. Moreover, if $Q$ is convex, then $1/3$ is the best possible.
\end{thm}

In the following, we obtain an multidimensional version of Theorem \ref{P9-thm-B} (b) for the operator valued analytic functions in the complete circular domain $Q$.
\begin{thm} \label{P9-thm-2.5}
If the series \eqref{him-vasu-P9-e-1.2} and \eqref{P9-e-2.5-adjust} converge in the domain $Q$ such that $\norm{f(z)} <1$ for all $z \in Q$ and $f(0)=a_{0}I$, $a_{0} \in \mathbb{C}$, $|a_{0}|<1$. If $G$ is given by \eqref{e-1.8}, then 
	\begin{equation} \label{P9-e-2.7}
	\sum_{k=0}^{\infty} \norm{P_{k}(z)} + G\left(\sum_{k=1}^{\infty}k \norm{P_{k}(z)}^2\right) \leq 1
	\end{equation}
	holds in the homothetic domain $(1/3)Q$, where the coefficients of $G$ satisfy
		\begin{equation} \label{P9-lem-3.20}
	8c_1\left(\frac{3}{8}\right)^2+24c_2\left(\frac{3}{8}\right)^4+\cdots+8(2l-1)c_l\left(\frac{3}{8}\right)^{2l} \leq 1.
	\end{equation}
	Moreover, if $Q$ is convex, then $1/3$ cannot be replaced by a larger quantity.
\end{thm}
\begin{cor} Let $f$ be as in Theorem \ref{P9-thm-2.5}. If $G$ is a monomial of one degree, then
\begin{equation}
\sum_{k=0}^{\infty} \norm{P_{k}(z)} + \frac{8}{9}\sum_{k=1}^{\infty}k \norm{P_{k}(z)}^2 \leq 1
\end{equation}
holds in the homothetic domain $(1/3)Q$. Moreover, if $Q$ is convex, then $1/3$ cannot be replaced by a larger quantity. 
\end{cor}
\section{Key Lemmas and their proofs}
We first obtain the upper and lower bounds of $m(\chi,r) $ for functions in $\mathcal{B}(\mathbb{D},\mathcal{B}(\mathcal{H}))$.
\begin{lem} \label{P9-lem-2.1}
	Let $f \in \mathcal{B}(\mathbb{D},\mathcal{B}(\mathcal{H}))$. Then
	\begin{enumerate}
		\item [\rm{(a)}] $m(\chi,r) \leq 1/\sqrt{1-r^2}$ for $r \in [0,1)$,
		\item [\rm{(b)}] $m(\chi,r) \geq (3-\sqrt{8(1-r^2)})/(1-r)$ for $r \in [1/3,1)$. 
	\end{enumerate}	
\end{lem}
\begin{pf} Let $f(z)=\sum_{k=0}^{\infty} A_{k}z^k$, where $A_{k} \in \mathcal{B}(\mathcal{H})$ for $k\in \mathbb{N}\cup \{0\}$.
\begin{enumerate}
	 \item [\rm{(a)}] In view of the Cauchy-Schwarz inequality, we have 
	\begin{equation} \label{P9-e-2.2}
	M_{r}(f) = \sum_{k=0}^{\infty} \norm{A_{k}} r^k \leq \sqrt{\sum_{k=0}^{\infty} \norm{A_{k}}^2} \, \sqrt{\sum_{k=0}^{\infty} r^{2k}} \,\,\,\,\, \,\,\,\mbox{for} \,\,\, r\in [0,1).
	\end{equation}
	From the given assumption that $f$ is in the unit ball of $H^{\infty}(\mathbb{D},\mathcal{B}(\mathcal{H}))$. That is, $\norm{f}_{H^{\infty}(\mathbb{D},\mathcal{B}(\mathcal{H}))} \leq 1$. In particular, we have $\norm{f}^2_{H^{2}(\mathbb{D},\mathcal{B}(\mathcal{H}))} =\sum_{k=0}^{\infty} \norm{A_{k}}^2\leq 1$.
Then the desired inequality follows from \eqref{P9-e-2.2}.\\
	\item [\rm{(b)}] To obtain the lower bound of $m(\chi,r)$, we consider the function $\psi_{a}$ which is defined by \eqref{P9-e-1.7}. 
	It is easy to see that
	$M_{r}(f) = (a+(1-2a^2)r)/(1-ar):=\pi(a)$ and hence, $m(\chi,r)\geq \pi(a)$ for $r\in [0,1)$. When $r\in [0,1/3]$, it is known that $m(\chi,r)\geq 1$. Now to obtain a better bound of $m(\chi,r)$, we shall maximize $\pi(a)$. For $r \geq 1/3$, the maximum value of $\pi(a)$ occurs at  $a=\left(2-\sqrt{2(1-r^2)}\right)/2r$. By substituting this value of $a$ in $\pi(a)$, we obtain $m(\chi,r) \geq (3-\sqrt{8(1-r^2)})/(1-r)$ for $r \in [1/3,1)$.
\end{enumerate}
\end{pf}

In the following result, we establish an operator valued analogues of Theorem \ref{P9-thm-B} (c) with the term $|f(z)|$ replaced by $\norm{f(z)}^p$ for $p\in (0,1]$ for functions $f$ in $ \mathcal{B}(\mathbb{D},\mathcal{B}(\mathcal{H}))$.
\begin{lem} \label{P9-lem-2.2}
	Let $f : \mathbb{D}  \rightarrow \mathcal{B}(\mathcal{H})$ be an operator valued bounded holomorphic function with the expansion $f(z)= \sum_{k=0}^{\infty} A_{k}z^k$ in $\mathbb{D}$ such that $A_{k} \in \mathcal{B}(\mathcal{H})$ for all $k \in \mathbb{N} \cup \{0\}$ and $A_{0}=a_{0}I$, $a_{0} \in \mathbb{C}$, $|a_{0}|<1$. If $\norm{f(z)} \leq 1$ in $\mathbb{D}$, then for $p \in (0,1]$, we have
	\begin{equation} \label{him-vasu-P9-e-2.5}
	\norm{f(z)}^p + \sum_{k=N}^{\infty} \norm{A_{k}}r^k \leq 1 \,\,\,\,\,\,\mbox{for} \,\,\,\,\,|z|=r \leq R_{N,p},
	\end{equation}
	where $R_{N,p}$ is the unique root in $(0,1)$ of the equation \eqref{P9-e-2.5}. The constant $R_{N,p}$ is the best possible.
\end{lem}
\begin{pf}
	Let $\mathcal{B}(\mathcal{H})$ be the set of bounded linear operators on a complex Hilbert space $\mathcal{H}$. Then $\mathcal{B}(\mathcal{H})$ is a complex Banach space with respect to the norm 
	\begin{equation*}
	\norm{A}=\sup_{h\in \mathcal{H}\setminus\{0\}} \frac{\norm{Ah}}{\norm{h}}= \sup_{h\in \mathcal{H}, \norm{h}=1} \norm{Ah},
	\end{equation*}
	where $A \in \mathcal{B}(\mathcal{H})$. Let $X^*$ be the dual space of $\mathcal{B}(\mathcal{H})$. For $x \in \mathcal{B}(\mathcal{H}) \setminus \{0\}$, let 
	\begin{equation*}
	T(x)= \{l_{x} \in X^*: l_{x}(x)=\norm{x} \,\,\, \mbox{and} \,\,\, \norm{l_{x}}=1 \}.
	\end{equation*}
	By the Hahn-Banach theorem, we have $T(x) \neq \phi$. Fix $z \in \mathbb{D}\setminus\{0\}$ and let $\nu = z/|z| \in \partial \mathbb{D}$. Without loss of generality, we assume $f(z) \neq 0$. Let $\tau \in \mathcal{B}(\mathcal{H})$ be any fixed point such that $\norm{\tau}=1$. We define the following holomorphic function in $\mathbb{D}$ by 
	\begin{equation*}
	g(\rho)=l_{\tau}(f(\rho \nu)) \,\,\,\,\, \mbox{for} \,\,\,\, \rho \in \mathbb{D},
	\end{equation*}   
	where $l_{\tau} \in T(\tau)$. Clearly, $|g(\rho)| \leq 1$ in $\mathbb{D}$. By applying the Schwarz-Pick lemma (often referred to as Lindel\"{o}f's inequality) (see \cite{alkhaleefah-2019}) to the complex-valued function $g$, we obtain 
	\begin{equation} \label{him-vasu-P9-e-2.6}
	|g(\rho)| \leq \frac{|g(0)|+|\rho|}{1+|g(0)| |\rho|}, \,\,\,\,\, \rho \in \mathbb{D}.
	\end{equation}
	It is easy to see that 
	\begin{equation} \label{him-vasu-P9-e-2.7}
	|g(0)| = |l_{\tau}(f(0))| \leq \norm{l_{\tau}} \norm{f(0)}=\norm{f(0)}.
	\end{equation}
	Let $|g(0)|=t$ and $|\rho|=s$. Then $s \in [0,1)$ and hence, the right hand side term $(t+s)/(1+ts)$ of \eqref{him-vasu-P9-e-2.6} is an increasing function in the variable $t \in [0,\infty)$. Therefore, by \eqref{him-vasu-P9-e-2.6} and \eqref{him-vasu-P9-e-2.7}, we obtain 
	\begin{equation} \label{him-vasu-P9-e-2.8}
	|l_{\tau}(f(\rho \nu))| \leq \frac{\norm{f(0)} + |\rho|}{1+ \norm{f(0)}|\rho|}.
	\end{equation}
	By choosing $\rho =|z|$ and $\tau =f(z)/\norm{f(z)}$ in \eqref{him-vasu-P9-e-2.8}, we obtain
	\begin{equation} \label{him-vasu-P9-e-2.9}
	\norm{f(z)} \leq \frac{\norm{f(0)} + |z|}{1+\norm{f(0)} |z|} \,\,\,\,\,\mbox{for} \,\,\,\,\, z \in \mathbb{D}.
	\end{equation} 
	Let $f(z)= \sum_{k=0}^{\infty} A_{k}z^k$ be analytic in $\mathbb{D}$ such that $A_{k} \in \mathcal{B}(\mathcal{H})$ for all $k \in \mathbb{N} \cup \{0\}$ and $A_{0}=a_{0}I$, $a_{0} \in \mathbb{C}$ with $|a_{0}|<1$. Then, in view of \eqref{him-vasu-P9-e-1.3} and \eqref{him-vasu-P9-e-2.9}, we have
	\begin{equation} \label{him-vasu-P9-e-2.5-a}
	\norm{f(z)}^p + \sum_{k=N}^{\infty} \norm{A_{k}}|z|^k \leq \left(\frac{|z|+ |a_{0}|}{1+|z||a_{0}|}\right)^p + (1-|a_{0}|^2) \frac{|z|^N}{1-|z|}=1+\Upsilon_{N,p}(|a_{0}|),
	\end{equation}
	where 
	\begin{equation*}
	\Upsilon_{N,p}(|a_{0}|)=  \left(\frac{r+|a_{0}|}{1+r|a_{0}|}\right)^p -\frac{1-r-(1-|a_{0}|^2)r^N}{1-r} \,\,\,\,\, \mbox{for} \,\,\,\, |z|=r
	\end{equation*}
	and $|a_{0}| \in [0,1]$. Let $|a_{0}|=\alpha$. To prove the inequality \eqref{him-vasu-P9-e-2.5}, it is enough to show that $\Upsilon_{N,p}(\alpha) \leq 0$ for $\alpha \in [0,1]$. That is, we have to show that $\Phi_{N,p}(\alpha) \leq 0$, where
	\begin{equation*}
	\Phi_{N,p}(\alpha)= (1-r)(r+\alpha)^p - (1+r\alpha)^p \left(1-r-(1-\alpha^2)r^N\right).
	\end{equation*}
	Clearly, $\Phi_{N,p}(1)=0$. Therefore, if we show that $\Phi_{N,p}$ is an increasing function in $\alpha$ under the condition \eqref{P9-e-2.5}, then we are done. A simple computation shows that 
	\begin{equation*}
	\Phi'_{N,p}(\alpha) = p(1-r)(r+\alpha)^{p-1} - pr\left(1-r-(1-\alpha^2)r^N\right)(1+r\alpha)^{p-1} - 2 \alpha r^N (1+r \alpha)^p
	\end{equation*}
	and
	\begin{align*}
	\Phi''_{N,p}(\alpha) &= p(p-1)(1-r)\left((r+\alpha)^{p-2} - r^2 (1+r\alpha)^{p-2} \right) + p(p-1)r^2 (1-\alpha^2)r^N (1+r\alpha)^{p-2} \\
	& \,\,\, \,\,\,\,- 2\alpha pr^{N+1}(1+r\alpha)^{p-1} - 2r^N (1+r \alpha)^p.
	\end{align*}
	Note that $(r+\alpha)^{p-2} - r^2 (1+r\alpha)^{p-2}$ is positive for $\alpha \in [0,1]$ and $p \in (0,1]$, which shows that the first term in the expression of $\Phi''_{N,p}(\alpha)$ is negative for $p \in (0,1]$. Since the other terms in the expression of $\Phi''_{N,p}(\alpha)$ are also negative for $p \in (0,1]$, it follows that $\Phi''_{N,p}(\alpha) \leq 0$ for all $\alpha \in [0,1]$ and $p \in (0,1]$. Thus, $\Phi'_{N,p}$ is a monotonically decreasing function in $\alpha \in [0,1]$ for $p \in (0,1]$, which gives that 
	\begin{equation}
	\Phi'_{N,p}(\alpha) \geq \Phi'_{N,p}(1)=p(1-r)^2 - 2r^N(1+r):= \Psi(r).
	\end{equation}
	We observe that $\Psi(r) \geq 0$ for $r \leq R_{N,p}$, where $R_{N,p}$ is the unique root of $\Psi(r)=0$. Indeed, $\Psi'(r)=-2p(1-r)-2\left(Nr^{N-1} + (N+1)r^N\right) \leq 0$ implies that $\Psi$ is a decreasing function in $r \in [0,1]$. On the other hand, since $\Psi (0)=p>0$ and $\Psi(1)=-4<0$, $\Psi$ has the unique root in $(0,1)$ and let $R_{N,p}$ be that root. Since $\Psi$ is a decreasing function in $r$, we have $\Psi(r) \geq \Psi(R_{N,p})=0$ for $r \leq R_{N,p}$. Therefore, $\Phi'_{N,p}(\alpha) \geq 0$ for $r \leq R_{N,p}$ and $p \in (0,1]$. Thus, $\Phi_{N,p}$ is an increasing function in $\alpha$ whenever $p \in (0,1]$. Therefore, $\Phi_{N,p}(\alpha) \leq \Phi_{N,p}(1)=0$ for $r \leq R_{N,p}$ for $p \in (0,1]$. This shows that $\Upsilon_{N,p}(\alpha) \leq 0$ for $r \leq R_{N,p}$, $p \in (0,1]$ and hence, the inequality \eqref{him-vasu-P9-e-2.5} follows from \eqref{him-vasu-P9-e-2.5-a}. 
	\vspace{2mm}
	
	To prove the sharpness of the radius $R_{N,p}$, we consider the function $\psi_{a}$ with $\psi_{a}(z)=\sum_{k=0}^{\infty} A_{k} z^k $ given by \eqref{P9-e-1.7}.
	A simple computation shows that 
	\begin{align} \label{P9-e-2.10}
	\norm{\psi_{a}(-r)}^p + \sum_{k=N}^{\infty} \norm{A_{k}}r^k &= \left(\frac{r+a}{1+ra}\right)^p + (1-a^2) \frac{a^{N-1}r^N}{1-ar} \\ \nonumber
	& = 1- \frac{(1-a)B_{N,p}(a,r)}{(1+ar)^p (1-ar)},
	\end{align}
	where
	\begin{equation*}
	B_{N,p}(a,r)= (1-ar)(1+ar)^p \left(\frac{1}{1-a}\left(1-\left(\frac{r+a}{1+ra}\right)^p\right) - \left(\frac{1+a}{1-ar}\right)a^{N-1}r^N\right).
	\end{equation*}
	We note that the right hand expression of \eqref{P9-e-2.10} is greater than or equals to $1$ if, and only if, $B_{N,p}(a,r) \leq 0$. By letting $a \rightarrow 1^{-}$, we obtain
	\begin{equation*}
	\lim\limits_{a \rightarrow 1^{-}} B_{N,p}(a,r) = (1-r)(1+r)^p \left(p\left(\frac{1-r}{1+r}\right) - \frac{2r^N}{1-r}\right)<0 
	\end{equation*}
	for $r>R_{N,p}$. This proves the sharpness of the constant $R_{N,p}$. This completes the proof.
\end{pf}
In the next result, we prove an operator valued analogues of Theorem \ref{P9-thm-B} (a).
\begin{lem} \label{P9-lem-2.13}
	Let $f : \mathbb{D}  \rightarrow \mathcal{B}(\mathcal{H})$ be an operator valued bounded holomorphic function with the expansion $f(z)= \sum_{k=0}^{\infty} A_{k}z^k$ in $\mathbb{D}$ such that $A_{k} \in \mathcal{B}(\mathcal{H})$ for all $k \in \mathbb{N} \cup \{0\}$ and $A_{0}=a_{0}I$, $a_{0} \in \mathbb{C}$ with $|a_{0}|<1$. If $\norm{f(z)} \leq 1$ in $\mathbb{D}$, then we have
	\begin{equation} \label{him-vasu-P9-e-2.13}
	\sum_{k=0}^{\infty} \norm{A_{k}}r^k + \left(\frac{1}{1+\norm{f(0)}} + \frac{r}{1-r}\right) \sum_{k=1}^{\infty} \norm{A_{k}}^2 r^{2k} \leq 1 \,\,\,\,\,\,\mbox{for} \,\,\,\,\,|z|=r\leq \frac{1}{3}.
	\end{equation}
	The constant $1/3$ is the best possible.
\end{lem}
\begin{pf}
	Let $f(z)= \sum_{n=0}^{\infty} A_{n}z^n$ in $\mathbb{D}$ with $\norm{f(z)} \leq 1$ in $\mathbb{D}$ such that $A_{n} \in \mathcal{B}(\mathcal{H})$ for all $n \in \mathbb{N} \cup \{0\}$ and $A_{0}=a_{0}I$, $a_{0} \in \mathbb{C}$. Set $\norm{A_{0}}=|a_{0}|=b \in [0,1]$. Then by \eqref{him-vasu-P9-e-1.3}, we obtain
	\begin{align*}&
	\sum_{k=0}^{\infty} \norm{A_{k}} r^k + \left(\frac{1}{1+|a_{0}|} + \frac{r}{1-r}\right) \sum_{k=1}^{\infty} \norm{A_{k}}^2 r^{2k}\\
	& \leq b+ (1-b^2)\left(\frac{r}{1-r}\right) + \left(\frac{1}{1+b} + \frac{r}{1-r}\right) (1-b^2)^2 \, \frac{r^2}{1-r^2}:=\Psi_{2}(b).
	\end{align*}
	Clearly, $\Psi_{2}$ can be expressed as $\Psi_{2}(b)=b+ \alpha (1-b^2) +\beta (1-b)(1-b^2)+ \gamma (1-b^2)^2$ for $b \in [0,1]$, where
	\begin{equation*}
	\alpha =\frac{r}{1-r},\, \beta =\frac{r^2}{1-r^2}\,\, \mbox{and} \,\,\, \gamma=\frac{r^3}{(1-r)(1-r^2)}.
	\end{equation*}
	Clearly, $\alpha,\beta$ and $\gamma$ are non-negative. We note that 
	\begin{align*}&
	\Psi_{2}'(b)\,=1-2\alpha b +\beta(3b^2 -2b-1) + 4\gamma(b^3 -b), \\
	& \Psi_{2}''(b)= -2\alpha + 2 \beta (3b -1) + 4\gamma (3 b^2 -1) \,\, \mbox{and} \,\,\, \Psi_{2}'''(b)=6\beta b +24 b \gamma.
	\end{align*}
	Since $\beta$ and $\gamma$ are non-negative, we have $\Psi_{2}'''(b)>0$ for all $b \in [0,1]$. That is, $\Psi_{2}''$ is an increasing function of $b$, which implies that
	\begin{equation*}
	\Psi_{2}''(b) \leq \Psi_{2}''(1)= -2 \alpha + 4\beta + 8\gamma= \frac{2r}{(1-r)(1-r^2)}\, \tau (r), 
	\end{equation*}
	where $\tau (r)= 4r^2 + 2r (1-r) - (1-r^2)= (1+r)(3r -1)$. It is easy to see that $\tau (r) \leq 0$ for $r \leq 1/3$. Therefore, $\Psi_{2}''(b) \leq 0$ for $b \in [0,1]$, which shows that $\Psi_{2}'$ is a decreasing function in $b \in [0,1]$. Thus, for $r \leq 1/3$, we obtain 
	\begin{equation*}
	\Psi_{2}'(b) \geq \Psi_{2}'(1) = 1-2\alpha = \frac{1-3r}{1-r}.
	\end{equation*}
	Clearly, for $r \leq 1/3$, we have $\Psi_{2}'(1)\geq 0$ and hence $\Psi_{2}'(b) \geq 0$ in $[0,1]$. This implies that $\Psi_{2}$ is an increasing function in $[0,1]$ and hence, we obtain $\Psi_{2}(b) \leq \Psi_{2}(1)=1$ for $r \leq 1/3$. This proves the desired inequality \eqref{him-vasu-P9-e-2.13}. 
	\par
	For the sharpness of the constant $1/3$, we consider the function $\psi_{a}$ with $\psi_{a}(z)=\sum_{k=0}^{\infty} A_{k} z^k $ is given by \eqref{P9-e-1.7}. A simple computation shows that
	\begin{align*}&
	\norm{A_{0}}+\sum_{k=1}^{\infty} \norm{A_{k}} r^k + \left(\frac{1}{1+|a_{0}|} + \frac{r}{1-r}\right) \sum_{k=1}^{\infty} \norm{A_{k}}^2 r^{2k}\\
	& = a+ \frac{1-a^2}{a}\sum_{k=1}^{\infty} a^n r^n + \left(\frac{1}{1+a} + \frac{r}{1-r}\right)  \, \frac{1-b^2}{b}\sum_{k=1}^{\infty}a^{2k} r^{2k}:= 1+(1-a)u(r):=v(r),
	\end{align*}
	where 
	\begin{equation*}
	u(r)= \frac{(1+a)r}{1-ar} + \left(\frac{1}{1+a} + \frac{r}{1-r}\right) \frac{(1+a)(1-a^2)r^2}{1-a^2 r^2} -1.
	\end{equation*} 
	We note that $u$ is strictly increasing function in $(0,1)$. Hence, for $r>1/3$, we have $u(r)>u(1/3)$. By letting $a$ very close to $1$, we obtain 
	\begin{equation*}
	\lim \limits_{a \rightarrow 1^{-}}u(1/3)= \frac{2(1/3)}{1-(1/3)} -1=0.
	\end{equation*}
	Therefore, $u$ is strictly positive function for $r >1/3$, as $a$ is very close to $1$. Hence, $v(r)>1$ for $r >1/3$, which shows that the radius $1/3$ is the best possible. 
\end{pf}

 In the following result, we establish an operator valued analogues of Theorem \ref{P9-thm-B} (b).
\begin{lem} \label{P9-lem-3.19}
	Let $f$ be as in Lemma \ref{P9-lem-2.13}. If $\norm{f(z)} \leq 1$ in $\mathbb{D}$, then we have 
	\begin{equation} \label{P9-e-3.20}
	\sum_{k=0}^{\infty}\norm{A_{k}}r^k+G(S(z))\leq 1\;\;\mbox{for}\;\; r\leq \frac{1}{3},
	\end{equation} where the coefficients of $ G(w) $ are given by \eqref{e-1.8} satisfy \eqref{P9-lem-3.20}. Furthermore, the constant $ 1/3 $ cannot be replaced by a larger quantity. Here $S$ is given by \eqref{P9-e-2.5-adjust}.
\end{lem}
\begin{pf}
	Let $\norm{A_{0}}=|a_{0}|=b \in [0,1]$. Then using the inequality \eqref{him-vasu-P9-e-1.3}, we obtain
	\begin{equation} \label{P9-e-2.18}
	S(z) \leq (1-b^2)^2 \sum_{k=1}^{\infty} kr^{2k}= (1-b^2)^2 \frac{r^2}{(1-r^2)^2}.
	\end{equation}
	The inequality \eqref{P9-e-2.18} along with \eqref{him-vasu-P9-e-1.3} gives 
	\begin{align} \label{P9-e-2.19}
	\norm{A_{0}} + \sum_{k=1}^{\infty} \norm{A_{k}} r^k + U(S(z)) 
	&\leq 
	b+ (1-b^2)\frac{r}{1-r} + \sum_{m=1}^{l} c_{m} \left(\frac{(1-b^2)r}{1-r^2}\right)^{2m}=1+H(r),
	\end{align}
	where 
	\begin{equation} \label{P9-e-2.20}
	H(r)= (1-b^2)\frac{r}{1-r} + \sum_{m=1}^{l} c_{m} \left(\frac{(1-b^2)r}{1-r^2}\right)^{2m} - (1-b).
	\end{equation}
	It is easy to see that $H(r)$ is an increasing function and hence $H(r) \leq H(1/3)$ for $r\leq 1/3$. A simple computation shows that 
	\begin{align*}
	H(1/3)=\frac{1-b^2}{2}\left(1+2F_l(b)-\frac{2}{1+b}\right):=\frac{1-b^2}{2}J(b),
	\end{align*}
	where 
	\begin{equation*}
	F_l(b)=\sum_{m=1}^{l}c_{l} (1-b^2)^{2m-1}\left(\frac{3}{8}\right)^{2m} \,\,\, \mbox{and} \,\,\,
	J(b)=1+2F_l(b)-\frac{2}{1+b}.
	\end{equation*}
	To show that $ H(r)\leq 0 $, it is enough to show that $ J(b)\leq 0 $ for $ b\in [0,1] $. 
	Since $ b\in [0,1] $, a simple computation shows that
	\begin{align*}
	b(1+b)^2 \left(\frac{3}{8}\right)^2\leq 4\left(\frac{3}{8}\right)^2, \,\, b(1+b)^2(1-b^2)^2 \left(\frac{3}{8}\right)^4\hspace{0.2in}\leq 4\left(\frac{3}{8}\right)^4,\cdots,\\ b(1+b)^2(1-b^2)^{2m-2} \left(\frac{3}{8}\right)^{2l}\leq 4\left(\frac{3}{8}\right)^{2l}.
	\end{align*}
	It is easy to see that
	\begin{align*}
	J^{\prime}(b)&=\frac{2}{(1+b)^2}\bigg(1-2c_1b(1+b)^2 \left(\frac{3}{8}\right)^2-6c_2 b(1+b)^2(1-b^2)^2 \left(\frac{3}{8}\right)^4\\&\quad\quad-\cdots-2(2l-1)c_l b(1+b)^2(1-b^2)^{2l-2} \left(\frac{3}{8}\right)^{2l}\bigg)\\&\geq  \frac{2}{(1+b)^2}\left(1-\left(8c_1\left(\frac{3}{8}\right)^2+24c_2\left(\frac{3}{8}\right)^4+\cdots+8(2l-1)c_l\left(\frac{3}{8}\right)^{2l}\right)\right)\\&\geq 0,\;\;\;\;\; \text{if}\;\; 8c_1\left(\frac{3}{8}\right)^2+24c_2\left(\frac{3}{8}\right)^4+\cdots+8(2l-1)c_l\left(\frac{3}{8}\right)^{2l}\leq 1.
	\end{align*}
	Therefore, $ J(b) $ is an increasing function in $ [0,1] $ if \eqref{P9-lem-3.20} holds.
	Hence, $ J(b)\leq J(1)= 0 $ for all $ b\in [0,1] $, which gives the desired inequality \eqref{P9-e-3.20}. To show the sharpness of the constant $1/3$, we consider the function $\psi_{a}$ given by \eqref{P9-e-1.7} with $\psi_{a}(z)=\sum_{k=0}^{\infty} A_{k} z^k $.
		A simple computation using \eqref{him-vasu-P9-e-1.3} shows that
	\begin{align*} 
	\sum_{n=0}^{\infty}{|a_n|}r^n+U(S(z))&=a +(1-a^2)\frac{r}{1-ar}+\frac{c_1r^2(1-a^2)^2}{(1-a^2r^2)^2}+\cdots+\frac{c_mr^{2m}(1-a^2)^{2m}}{(1-a^2r^2)^{2m}} \\& :\;=1-(1-a)\Phi_1(r),
	\end{align*}
	where 
	\begin{align*} 
	\Phi_1(r)&=-\frac{(1+a)r}{1-ar}-\frac{c_1r^2(1-a)(1+a)^2}{(1-a^2r^2)^2}- \cdots -\frac{c_mr^{2m}(1-a)^{2m-1}(1+a)^{2m}}{(1-a^2r^2)^{2m}} +1.
	\end{align*}
	It is not difficult to show that $ \Phi_1(r) $ is strictly decreasing function of $ r$ in $(0,1) $. Therefore,  for $ r>1/3 $, we have $ \Phi_1(r)<\Phi_1(1/3) $. A simple computation shows that $\lim_{a\rightarrow 1}\Phi_1(1/3)= 0$.
	Therefore  $ \Phi_1(r)< 0 $ for $ r>1/3 $. Hence $ 1-(1-a)\Phi_1(r)>1 $ for $ r>1/3, $ which shows that $ 1/3 $ is the best possible. 
	This completes the proof.
\end{pf}
\section{Proofs of the main results}
\begin{pf} [{\bf Proof of Theorem \ref{him-vasu-P9-thm-2.1}}]
	To obtain the inequality \eqref{him-vasu-P9-e-2.1}, we convert the multidimensional power series \eqref{him-vasu-P9-e-1.2} into the power series of one complex variable and we want to make use of Lemma \ref{P9-lem-2.2}. Let $L=\{z=(z_{1}, \ldots,z_{n}):z_{j}=a_{j}t, j=1,2,\ldots,n; t \in \mathbb{C}\}$ be a complex line. Then, in each section of the domain $Q$ by the line $L$, the series \eqref{him-vasu-P9-e-1.2} turns into the following power series of complex variable $t$:
	\begin{equation} \label{P9-e-2.9}
	f(at)=\sum_{k=0}^{\infty} P_{k}(a) t^k= f(0) + \sum_{k=1}^{\infty} P_{k}(a) t^k.
	\end{equation}
	Since $\norm{f(at)}<1$ for $t \in \mathbb{D}$ and $f(0)=a_{0}I$, $a_{0}\in \mathbb{C}$, by making use of Lemma \ref{P9-lem-2.2}, we obtain
	\begin{equation} \label{P9-e-2.11}
	\sum_{k=0}^{\infty} \norm{P_{k}(a) t^k}<1
	\end{equation}
	for $z$ in the section $L \bigcap \left(\frac{1}{3}Q\right)$. Since $L$ is an arbitrary complex line passing through the origin, the inequality \eqref{P9-e-2.11} is just \eqref{him-vasu-P9-e-2.1}. For the sharpness of the constant $1/3$, let the domain $Q$ be convex. Then $Q$ is an intersection of half spaces 
	\begin{equation*}
	Q= \bigcap \limits _{a \in J} \{z=(z_{1}, \ldots, z_{n}): \real (a_{1}z_{1} + \cdots + a_{n}z_{n})<1\}
	\end{equation*}
	for some $J$. Because $Q$ is circular, we obtain 
	\begin{equation*}
	Q= \bigcap \limits _{a \in J} \left\{z=(z_{1}, \ldots, z_{n}): \abs {a_{1}z_{1} + \cdots + a_{n}z_{n}}<1\right\}.
	\end{equation*}
	Therefore, to show the constant $1/3$ is the best possible, it is enough to show that $1/3$ cannot be improved for each domain $G_{a}=\{z=(z_{1}, \ldots, z_{n}): \abs {a_{1}z_{1} + \cdots + a_{n}z_{n}}<1\}$. In view of Theorem \ref{him-vasu-P9-thm-1.1}, for some $b \in [0,1)$, there exists a function $\psi_{b}:\mathbb{D} \rightarrow \mathcal{B}(\mathcal{H})$ defined by \eqref{P9-e-1.7} with $\norm{\psi_{b}(\xi)}<1$ for $\xi \in \mathbb{D}$, but for any $|\xi|=r>1/3$, \eqref{him-vasu-P9-e-2.1} fails to hold in the disk $\mathbb{D}_{r}=\{\xi : |\xi|<r\}$. On the other hand, consider the function $\phi: G_{a} \rightarrow \mathbb{D}$ defined by $\phi(z)=a_{1}z_{1} + \cdots + a_{n}z_{n}$. Thus, the function $f(z)=(\psi_{b} \circ \phi) (z)$ gives the sharpness of the constant $1/3$ for each domain $G_{a}$. This completes the proof.  
\end{pf}
\vspace{-5mm}

\begin{pf} [{\bf Proof of Lemma \ref{P9-lem-2.7}}]
	In view of Lemma \ref{P9-lem-2.1} (a), with $r=1/\sqrt{2}$, we have $m(\chi,1/\sqrt{2}) \leq \sqrt{2}$ and by  Lemma \ref{P9-lem-2.1} (b) with $r=1/\sqrt{2}$, we obtain $m(\chi,1/\sqrt{2}) \geq \sqrt{2}$. These together give $(1/\sqrt{2}) \, m(\chi,1/\sqrt{2}) =1$. Hence by \eqref{P9-e-1.3}, we conclude that $R_{1}=1/\sqrt{2}$.
\end{pf}
\vspace{-5mm}

\begin{pf} [{\bf Proof of Theorem \ref{him-vasu-P9-thm-2.2}}]
	In view of Lemma \ref{P9-lem-2.7}, by going the similar lines of argument as in Theorem \ref{him-vasu-P9-thm-2.1}, we can easily show that $\sum_{k=1}^{\infty} \norm{P_{k}(z)} <1$ in the homothetic domain $(1/\sqrt{2})Q$. To show the constant $1/\sqrt{2}$ is the best possible, let the domain $Q$ be convex. Therefore, by the analogues proof of Theorem \ref{him-vasu-P9-thm-2.1}, it is enough to prove that the constant $1/\sqrt{2}$ cannot be improved in each domain $G_{a}=\{z=(z_{1}, \ldots, z_{n}): \abs {a_{1}z_{1} + \cdots + a_{n}z_{n}}<1\}$. Since $1/\sqrt{2}$ is the best possible in Lemma \ref{P9-lem-2.7}, there exists an analytic function $h :\mathbb{D} \rightarrow \mathcal{B}(\mathcal{H})$ defined by 
	\begin{equation*}
	h(\xi)= \xi \left(\frac{\frac{1}{\sqrt{2}}-\xi}{1-\frac{\xi}{\sqrt{2}}}\right)I= \sum_{k=1}^{\infty} A_{k} \xi ^k \,\,\,\,\, \mbox{for} \,\,\,\,\, \xi \in \mathbb{D},
	\end{equation*}
	where $A_{1}=1/\sqrt{2}$ and $A_{k}=-(1/2) (1/\sqrt{2})^{k-2}$ for $k \geq 2$ such that $\norm{h(\xi)}<1$ in $\mathbb{D}$ and $h(0)=0$. But for any $|\xi|=r>1/\sqrt{2}$,
	\begin{equation*}
	\sum_{k=1}^{\infty} \norm{A_{k}}r^k = \frac{r/\sqrt{2}}{1-(r/\sqrt{2})}>1,
	\end{equation*}
	which shows that \eqref{P9-lem-2.7} fails to hold in the disk $\mathbb{D}_{r}=\{z:|z|<r\}$. Therefore, the function $f(z)=(h \circ \phi) (z)$ gives the sharpness of the constant $1/\sqrt{2}$ in each domain $G_{a}$, where $\phi: G_{a} \rightarrow \mathbb{D}$ defined by $\phi(z)=a_{1}z_{1} + \cdots + a_{n}z_{n}$. This completes the proof.
\end{pf}
\vspace{-5mm}

\begin{pf}[{\bf Proof of Theorem \ref{P9-thm-2.3}}]
In view of Lemma \ref{P9-lem-2.2} and the analogues proof of Theorem \ref{him-vasu-P9-thm-2.1}, as well as from \eqref{him-vasu-P9-e-2.5}, we can easily obtain the inequality \eqref{P9-e-2.4} in the homothetic domain $(R_{N,p})Q$, where $R_{N,p}$ is the positive root in $(0,1)$ of the equation \eqref{P9-e-2.5}. To prove the constant $R_{N,p}$ is the best possible whenever $Q$ is convex, in view of the analogues proof of Theorem \ref{him-vasu-P9-thm-2.1}, it is enough to show that $R_{N,p}$ cannot be improved for each domain $G_{a}=\{z=(z_{1}, \ldots, z_{n}): \abs {a_{1}z_{1} + \cdots + a_{n}z_{n}}<1\}$. Since $R_{N,p}$ is the best possible in Lemma \ref{P9-lem-2.2}, there exists an analytic function $h:\mathbb{D} \rightarrow \mathcal{B}(\mathcal{H})$ such that $\norm{h(\xi)} <1$ in $\mathbb{D}$, but \eqref{P9-e-2.4} fails to hold in the disk $\mathbb{D}_{r}=\{z: |z|<r\}$ for each $|\xi|=r>R_{N,p}$. Thus, the function $f(z)=(h \circ \phi) (z)$ gives the sharpness of the constant $R_{N,p}$ in each domain $G_{a}$, where $\phi: G_{a} \rightarrow \mathbb{D}$ is defined by $\phi(z)=a_{1}z_{1} + \cdots + a_{n}z_{n}$.
\end{pf}
\vspace{-5mm}

\begin{pf}[{\bf Proof of Theorem \ref{P9-thm-2.4}}]
	In view of Lemma \ref{P9-lem-2.13} and the analogues proof of Theorem \ref{him-vasu-P9-thm-2.1}, as well as from \eqref{him-vasu-P9-e-2.13}, we can obtain the inequality \eqref{P9-e-2.5-a} in the homothetic domain $(1/3)Q$. To prove the constant $1/3$ is the best possible when $Q$ is convex, in view of the analogues proof of Theorem \ref{him-vasu-P9-thm-2.1}, it is enough to show that $1/3$ cannot be improved for each domain $G_{a}=\{z=(z_{1}, \ldots, z_{n}): \abs {a_{1}z_{1} + \cdots + a_{n}z_{n}}<1\}$. Thus, the function $f(z)=(h \circ \phi) (z)$ gives the sharpness of the constant $R_{N,p}$ in each domain $G_{a}$, where $\phi: G_{a} \rightarrow \mathbb{D}$ is defined by $\phi(z)=a_{1}z_{1} + \cdots + a_{n}z_{n}$.
\end{pf}
\vspace{-5mm}

\begin{pf}[{\bf Proof of Theorem \ref{P9-thm-2.5}}]
	Using Lemma \ref{P9-lem-3.19} and the analogues proof of Theorem \ref{him-vasu-P9-thm-2.1}, from \eqref{P9-e-3.20}, we can easily obtain the inequality \eqref{P9-e-2.7} in the homothetic domain $(1/3)Q$, where the coefficients of $G$ satisfy \eqref{P9-lem-3.20}. When $Q$ is convex, to prove the constant $1/3$ is the best possible, in view of the analogues proof of Theorem \ref{him-vasu-P9-thm-2.1}, it is enough to show that $1/3$ cannot be improved for each domain $G_{a}=\{z=(z_{1}, \ldots, z_{n}): \abs {a_{1}z_{1} + \cdots + a_{n}z_{n}}<1\}$. Thus, the function $f(z)=(h \circ \phi) (z)$ gives the sharpness of the constant $R_{N,p}$ in each domain $G_{a}$, where $\phi: G_{a} \rightarrow \mathbb{D}$ defined by $\phi(z)=a_{1}z_{1} + \cdots + a_{n}z_{n}$.
\end{pf}

\noindent\textbf{Acknowledgment:} 
The first author is supported by SERB-CRG and the second author is supported by CSIR (File No: 09/1059(0020)/2018-EMR-I), New Delhi, India.

\end{document}